\documentclass{gtmon_a}
\pdfoutput=1

\usepackage{amscd,amsxtra}

\proceedingstitle{Proceedings of the Nishida Fest (Kinosaki 2003)}
\conferencestart{28 July 2003}
\conferenceend{8 August 2003}
\conferencename{International Conference in Homotopy Theory}
\conferencelocation{Kinosaki, Japan}

\editor{Matthew Ando}
\givenname{Matthew}
\surname{Ando}

\editor{Norihiko Minami}
\givenname{Norihiko}
\surname{Minami}

\editor{Jack Morava}
\givenname{Jack}
\surname{Morava}

\editor{W Stephen Wilson}
\givenname{W Stephen}
\surname{Wilson}

\title{Interactions of strings and equivariant homology theories}

\author{Shingo Okuyama}
\givenname{Shingo}
\surname{Okuyama}
\address{Takuma National College of Technology\\\newline
Kagawa 769-1192\\Japan}
\email{okuyama@dc.takuma-ct.ac.jp}
\urladdr{}

\author{Kazuhisa Shimakawa}
\givenname{Kazuhisa}
\surname{Shimakawa}
\address{Takuma National College of Technology\\\newline
Kagawa 769-1192\\Japan}
\email{kazu@math.okayama-u.ac.jp}
\urladdr{}

\volumenumber{10}
\issuenumber{}
\publicationyear{2007}
\papernumber{19}
\startpage{333}
\endpage{346}

\doi{}
\MR{}
\Zbl{}

\arxivreference{}


\received{13 September 2004}
\revised{13 May 2006}
\accepted{}
\published{18 April 2007}
\publishedonline{18 April 2007}
\proposed{}
\seconded{}
\corresponding{}
\version{}

\makeatletter
\def\cnewtheorem#1[#2]#3{\newtheorem{#1}{#3}[section]
\expandafter\let\csname c@#1\endcsname\c@Theorem}

\AtBeginDocument{\let\bar\wbar\let\tilde\wtilde\let\hat\what}
\def\S{Section }
\let\tsty\textstyle

\newtheorem{Theorem}{Theorem}[section]
\cnewtheorem{Proposition}[Theorem]{Proposition}
\cnewtheorem{Lemma}[Theorem]{Lemma}
\cnewtheorem{Corollary}[Theorem]{Corollary}
\theoremstyle{definition}
\cnewtheorem{Definition}[Theorem]{Definition}
\cnewtheorem{Example}[Theorem]{Example}
\theoremstyle{remark}
\cnewtheorem{Remark}[Theorem]{Remark}

\makeatother  

\makeautorefname{Theorem}{Theorem}
\makeautorefname{Proposition}{Proposition}
\makeautorefname{Lemma}{Lemma}

\numberwithin{equation}{section}

\newcommand{\propref}[1]{\fullref{#1}}

\DeclareMathOperator{\Map}{Map}
\DeclareMathOperator{\Top}{Top}

\newcommand{\graded}[2][\bullet]{{#2}_{#1}}
\newcommand{\thickr}[1]{\|{#1}\|}

\newcommand{\Rinfty}{\mathbb{R}^{\infty}}
\newcommand{\cQ}[1]{\mathcal{Q}(#1)}
\newcommand{\RS}[1]{\lvert{S_{\bullet}{#1}}\rvert}
\newcommand{\BD}{\mbox{\it BD\/}}
\newcommand{\Loop}{\Omega}
\newcommand{\simp}[1]{{#1}_{\bullet}}

\begin{document}

\begin{abstract}
We introduce the notion of the space of parallel strings with partially
summable labels, which can be viewed as a geometrically constructed group
completion of the space of particles with labels.  We utilize this to
construct a machinery which produces equivariant generalized homology
theories from such simple and abundant data as partial monoids.
\end{abstract}

\maketitle

%
\section{Introduction}

In \cite{config} we attached to any pair of a Euclidean space $V$ and
a partial abelian monoid $M$ a space $C(V,M)$ whose points are pairs
consisting of a finite subset $c$ of $V$ and a map $a \co  c \to M$,
but $(c,a)$ is identified with $(c',a')$ if $c \subset c'$, $a'|c =
a$, and $a'(v) = 0$ for $v \not\in c$.
Any such pair $(c,a)$ can be identified with the set consisting of
``labeled particles'' $(v,a(v))$, $v \in c$. 
Suppose $V$ is an orthogonal $G$--module for some finite group $G$ and
$M$ admits a $G$--action compatible with partial sum operations.  Then
$C(V,M)$ is a $G$--space with respect to the $G$--action
\[
g(c,a) = (gc,gag^{-1}), \quad g \in G,\ (c,a) \in C(V,M).
\]
Let $I(\R)$ be the space of finite disjoint unions of bounded
intervals in the real line.
Then $I(\R)$ is a
partial abelian monoid with partial sum operation given by
superimposition.
Let us denote $I(V,M) = C(V,I(\R) \wedge M)$ for any partial abelian
monoid with $G$--action $M$.  Observe that under the correspondence
\[
\textstyle%
a \co  c \to I(\R) \quad \mapsto \quad \bigcup_{v \in c}\{v\} \times
a(v) \subset V \times \R
\]
any map from a finite subset of $V$ to $I(\R)$ can be identified with
a finite disjoint union of bounded subsets of the form $\{v\} \times J
\subset V \times \R$, where $J$ is a bounded interval.  We call such
$\{v\} \times J$ a {\em string\/} over $v$.  Thus $I(V,M)$ can be
regarded as the space consisting of finite sets of pairwise disjoint
labeled strings whose members over the same point in $V$ has the same
label in $M$.

The aim of this paper is to show that if $V$ is sufficiently large
then there is a $G$--equivariant group completion map $C(V,M) \to
I(V,M)$ and also that the correspondence $X \mapsto \pi_n I(V,X \wedge
M)$, $n \geq 0$, extends to an $RO(G)$--graded generalized homology
theory.

To state the precise results, let $\Top(G)$ be the category of all
pointed $G$--spaces and all pointed maps with $G$ acting on maps by
conjugation.
In \cite{infinite} we have shown that any $G$--equivariant continuous
functor $T \co  \Top(G) \to \Top(G)$ such that $T(*) = *$ is
associated with pairings
\(
X \wedge TY \to T(X \wedge Y),\ TX \wedge Y \to T(X \wedge Y)
\)
natural in both $X$ and $Y$.
Therefore, $T$ preserves $G$--homotopies and there is a natural
transformation $S^W \wedge T(X) \to T(S^W \wedge X)$ for any
orthogonal $G$--module $W$, where $S^W$ is the one point
compactification of $W$.

Suppose $V$ is linearly and equivariantly isometric to the direct
product of countably many copies of the regular representation of $G$
over the real number fields.  Such a $G$--module $V$ is called a
$G$--universe.  Now the main results can be stated as follows.

\begin{Theorem}\label{thm:group-completion}
  There is a diagram consisting of maps of Hopf $G$--spaces
  \[
  C(V,M) \xleftarrow{\lambda} I_+(V,M) \xrightarrow{\rho} I(V,M)
  \]
  satisfying the conditions below.
  \begin{enumerate}
  \item $\lambda$ is a $G$--homotopy equivalence.
  \item $\rho$ is an equivariant group completion, that is to say, it
    restricts to a group completion map $I_+(V,M)^H \to I(V,M)^H$ for
    every subgroup $H$ of $G$.
  \end{enumerate}
\end{Theorem}

\begin{Theorem}\label{thm:infinite-delooping}
  The correspondence $X \to I(V,X \wedge M)$ is a $G$--equivariant
  continuous functor of $\Top(G)$ into itself and we have the following:
  \begin{enumerate}
  \item For any orthogonal $G$--module $W$ the natural map
    \[ I(V,X \wedge M) \to \Omega^W I(V,\Sigma^W X \wedge M) \]%
    adjoint to $S^V \wedge I(V,X \wedge M) \to I(V,S^W \wedge X \wedge
    M)$ is a weak $G$--equivalence.
  \item There exists an $RO(G)$--graded homology theory
    $\graded{h}^G(-)$ such that
    \[
    h^G_n(X) = \pi_n I(V,X \wedge M)^G
    \]
    holds for any $X$ and $n \geq 0$.
  \end{enumerate}
\end{Theorem}

These theorems enable us to construct equivariant generalizations of
several popular homology theories.  For example, consider the simplest
case $M = S^0$.  Then $C(V,X)$ is the usual configuration space, and
hence its group completion $I(V,X)$ is weakly $G$--equivalent to the
equivariant infinite loop space $\Omega^{V}\Sigma^{V}X$ by
\cite[Theorem (1.18)]{caruso-waner:approximation}. Thus we obtain the
$G$--equivariant stable homotopy theory in this case.
On the other hand, if we take arbitrary positive integers as labels
then we obtain an $RO(G)$--graded homology theory extending the
ordinary homology $\widetilde{H}_n(X/G,\mathbb{Z})$.  (Compare
Lewis, May and McClure \cite{lewis-may-mcclure}.)
$K$--theory type examples also occur from our method, which will be
discussed in a future paper.
\section[Partial abelian monoids with G-action]{Partial abelian monoids with $G$--action}

\begin{Definition}
  \label{def:pam}
  A pointed $G$--space $M$ is called a partial abelian monoid with
  $G$--action, or {\em $G$--partial monoid\/} for short, if for every $n
  \geq 0$ there are $G$--invariant subsets $M_n$ of $M^n$ and
  $G$--maps  
  \begin{center}
    $M_n \to M,\quad (a_1,\dots,a_n) \mapsto a_1 + \cdots + a_n$
  \end{center}
  satisfying the conditions below.
  \begin{enumerate}
  \item $M_0 \to M$ is the inclusion of the basepoint $0$ of $M$.
  \item $M_1 \to M$ is the identity of $M$.
  \item Let $J_1,\ \cdots,\ J_r$ be disjoint subsets of $\{1,\dots,n\}$
    such that $J_1 \cup\cdots\cup J_r = \{1,\dots,n\}$, and let
    $(a_1,\dots,a_n)$ be an element of $M^n$ such that $(a_j)_{j \in
      J_k}$ belongs to $M_{n_k}$, where $n_k$ is the cardinality of
    $J_k$.
    Then $(a_1,\dots,a_n) \in M_n$ if and only if $\big(\tsty\sum_{j \in
        J_1} a_j,\dots,\tsty\sum_{j \in J_r} a_j \big) \in M_r$, and we
    have
    \[
    \textstyle a_1 + \cdots + a_n = \sum_{j \in J_1} a_j + \cdots +
    \sum_{j \in J_r} a_j
    \]
    if either side of the equation makes sense.
  \end{enumerate}   
\end{Definition}

Among the examples we have the following:
\begin{enumerate}
\item Let $M$ be a $G$--invariant subset of a topological abelian
  group on which $G$ acts through group homomorphisms.  Suppose $M$
  contains the unit $0$.  Then $M$ is a $G$--partial monoid with
  respect to the subsets
  \[
  M_n = \{ (a_1,\dots,a_n) \in M^n \mid a_1 + \cdots +a_n \in M \}.
  \]
  More generally, any $G$--invariant subset of a $G$--partial monoid
  that contains $0$ is again a $G$--partial monoid.

\item Any pointed $G$--space $X$ is a $G$--partial monoid with respect
  to folding maps $X_n = X
  \vee \cdots \vee X \to X$.  In fact, this is a special case of the
  previous example, as $X$ is a $G$--invariant subset of the infinite
  symmetric product $\mbox{SP}^{\infty}X$.
  
\item Let $V$ be an infinite dimensional real inner product space on
  which $G$ acts through linear isometries.  Then the Grassmannian
  $\mbox{Gr}(V)$ of finite-dimensional subspaces of $V$ is a
  $G$--partial monoid with respect to the inner direct sum operation
  $\mbox{Gr}(V)_n \to \mbox{Gr}(V)$, where $\mbox{Gr}(V)_n$ is defined
  to be the subset consisting of those $(W_1,\dots,W_n)$ such that
  $W_i \perp W_j$ if $i \neq j$.
\end{enumerate}

\begin{Definition}
  For given $G$--partial monoids $M$ and $N$, their smash product $M
  \wedge N$ is a $G$--partial monoid such that $(M \wedge N)_n$ is
  the subset consisting of those $n$--tuples that can be summed up to
  an element of $M \wedge N$ by using the distributivity relations:
  \[\begin{array}{rcll}
    c_1 \wedge d + \cdots + c_k \wedge d &=& (c_1 + \cdots + c_k)
    \wedge d, &\quad(c_1,\dots,c_k) \in M_k
    \\
    c \wedge d_1 + \cdots + c \wedge d_l &=& c \wedge (d_1 + \cdots
    + d_l), &\quad(d_1,\dots,d_l) \in N_l
  \end{array}\]
\end{Definition}

\begin{Example}
  If $X$ is a pointed $G$--space and $M$ is a $G$--partial monoid,
  then $X \wedge M$ is a $G$--partial monoid such that
  \[
  (X \wedge M)_n = X \wedge M_n
  \]
  holds for every $n \geq 0$.
\end{Example}

For any orthogonal $G$--module $V$, the labeled configuration space
$C(V,M)$ is a $G$--partial monoid with respect to the partial sum
operations
\[
\textstyle C(V,M)_n \to C(V,M), \quad ((c_1,a_1),\dots,(c_n,a_n))
\mapsto (\bigcup c_i,\bigcup a_i)\,.
\]
Here $C(V,M)_n$ consists of those $n$--tuples $((c_i,a_i)) \in
C(V,M)^n$ such that for every $x \in \bigcup c_i$ the sum $\sum_{i \in
  \Lambda(x)}a_i(x)$ exists, where $\Lambda(x) = \{i \mid x \in c_i\}$,
and $\bigcup a_i$ denotes the map $x \mapsto \sum_{i \in
  \Lambda(x)}a_i(x)$.
Moreover, if $V$ is a $G$--universe then $C(V,M)$ is a homotopy
associative and homotopy commutative Hopf $G$--space.  To see this, let
us consider the functor
\[
P \mapsto A(P) = C(V,P \wedge M)
\]
from finite pointed sets to pointed $G$--spaces.
For each $p \in P$, let $\delta_p$ be the pointed map
$P \to \mathbf{1} = \{0,1\}$ such that
$\delta_p^{-1}(1) = \{p\}$ if $p$ is not the basepoint of $P$,
and let $\delta_p$ be the constant map if $p$ is the basepoint.
Then the $G$--map
\[
\delta \co  A(P) \to \Map_0(P,A(\mathbf{1})),
\quad a \mapsto (p \mapsto A(\delta_p)(a))
\]
has a $G$--homotopy inverse $\psi \co  \Map_0(P,A(\mathbf{1})) \to
A(P)$ defined as follows.

Since $V$ is a $G$--universe, there exist an embedding of $P - \{0\}$
into $V^G$ and a $G$--linear isometry $V \times V \to V$.  Hence we can
construct a $G$--equivariant embedding of $(P - \{0\}) \times V$ into
$V$.
For any $f \in \Map_0(P,A(\mathbf{1}))$ let us write $f(p) =
(c(p),a(p))$ and put $\psi(f) = (\hat{c},\hat{a}) \in A(P)$,
where $\hat{c}$ is the image of
\[
\textstyle {\bigcup}_{p \in P - \{0\}}\ \{p\} \times c(p)
\]
under the $G$--equivariant embedding $(P - \{0\}) \times V \to V$ and
$\hat{a} \co  \hat{c} \to P \wedge M$ is induced by the composite
maps
\[
c(p) \xrightarrow{a(p)} M = \mathbf{1} \wedge M \xrightarrow{\iota_p
  \wedge 1} P \wedge M
\]
where $\iota_p$ is a pointed map $\mathbf{1} \to P$ such that
$\iota_p(1) = p$.

Therefore, $A$ is a $G$--equivariant $\Gamma$--space in the sense of
Segal.  Hence the following proposition holds.

\begin{Proposition}
  $C(V,M)$ is a homotopy associative and homotopy commutative Hopf
  $G$--space with unit $\emptyset \in C(V,M)^G$.
\end{Proposition}

Note that Hopf $G$--space multiplication $\mu$ of $C(V,M)$ is given by
the composite
\[
C(V,M)^2 \xrightarrow[\simeq]{\psi} C(V,M \vee M)
\xrightarrow{\nabla_*} C(V,M)
\]
where $\nabla_*$ is induced by the folding map $M \vee M \to M$.

\begin{Definition}
  \label{def:homotopically-invertible}
  A $G$--partial monoid $M$ is {\em homotopically invertible\/} if there
  exist a map of $G$--partial monoids $\tau \co  M \to M$, called a
  {\em homotopy inversion\/}, and a $G$--homotopy $h_t \co  M \to M^2$
  ($0 \leq t \leq 1$) satisfying the conditions below.
  \begin{enumerate}
  \item For every $t \in [0,1]$, $h_t$ is a map of $G$--partial monoid.
  \item $h_0 = (1,\tau)$, ie we have $h_0(a) = (a,\tau(a))$ for any
    $a \in M$.
  \item $h_1$ factors through a map $h'_1 \co  M \to M_2$ and the
  composite $\smash{M \xrightarrow{\smash{h'_1}\,} M_2 \xrightarrow{\Sigma} M}$ is
  $G$--homotopic through maps of $G$--partial monoids to the
    constant map.
  \end{enumerate}
\end{Definition}

\begin{Proposition}
  If $V$ is a $G$--universe and if $M$ is homotopically invertible then
  $C(V,M)$ is a grouplike Hopf $G$--space.
\end{Proposition}

\begin{proof}
  Let $\tau_* \co  C(V,M) \to C(V,M)$ be the map induced by the
  homotopy inversion of $M$.  To see that $C(V,M)$ is grouplike, it
  suffices to show that the composite
  \[
  C(V,M) \xrightarrow{(1,\tau_*)} C(V,M)^2 \xrightarrow{\mu} C(V,M)
  \]
  is $G$--homotopic to the constant map with value $\emptyset$.
  Let us regard $M \times M$ as a $G$--partial monoid such that
  $(M \times M)_n = M_n \times M_n$ for $n \geq 0$.  Then we have a
  diagram of pointed $G$--spaces
  \begin{equation*}
    \begin{CD}
      C(V,M) @>{(1,\tau)_*}>> C(V,M^2) @>{(p_{1*},p_{2*})}>> C(V,M)^2
      \\
      @V{h'_{1*}}VV @AAA @A{\simeq}A{\delta}A
      \\
      C(V,M_2) @= C(V,M_2) @<<< C(V,M \vee M)
      \\
      @. @V{\Sigma_*}VV @VV{\nabla_*}V
      \\
      @. C(V,M) @= C(V,M)
    \end{CD}
  \end{equation*}
  in which $p_{1*}$ and $p_{2*}$ are induced by the projections $M^2
  \to M$ onto the first and the second factors, respectively, and
  unnamed arrows are induced by the inclusions of $G$--partial monoids.
  Clearly, the right hand side squares are commutative, and the upper
  left square commutes up to $G$--homotopy.
  Since $\delta$ has a $G$--homotopy inverse $\psi$ and since $\psi$
  restricts to a $G$--homotopy inverse to the map $C(V,M \vee M) \to
  C(V,M^2)$ induced by the inclusion $M \vee M \to M^2$, all the
  arrows constituting the upper right square are $G$--homotopy
  equivalences.  Thus we have
  \[
  \mu(1,\tau_*) = \nabla_*\psi(p_{1*},p_{2*})(1,\tau)_* \simeq
  \Sigma_*h'_{1*} \simeq \emptyset.\proved
  \]
\end{proof}
\section{The space of strings with labels}

As usual, the symbols $[a,b]$, $[a,b)$, $(a,b]$, $(a,b)$ represent
bounded intervals in the real line, and $b - a$ is called the length
of the interval.
The space of intervals $I(\R)$ consists of those unions $P = J_1 \cup
\cdots \cup J_r$ of finite number of pairwise disjoint bounded
intervals.
It is topologized in such a way that such operations as isotopy moves,
concatenation of two disjoint intervals that have a connected union
(eg\ $[a,c) \cup [c,b] = [a,b]$), and deletion of a half-open
interval when its length tends to $0$ are all continuous.
Let $I(\R)_n$ be the subset of $I(\R)^n$ consisting of those
$n$--tuples $(P_1,\dots,P_n)$ that are pairwise disjoint.  Then $I(\R)$
is a partial abelian monoid with respect to these $I(\R)_n$ and maps
\[
I(\R)_n \to I(\R), \quad (P_1,\dots,P_n) \mapsto P_1 \cup\cdots\cup
P_n.
\]
Details are given in Okuyama \cite{okuyama:intervals}, where $I(\R)$ is denoted by
$I_1(S^0)$.

\begin{Lemma}
  \label{lem:IR-is-homotopically-invertible}
  $I(\R)$ is a homotopically invertible partial abelian monoid.
\end{Lemma}

\begin{proof}
  Given a bounded interval $J$, let $\tau{J}$ denote the complement of
  the boundary of $-J$ in its closure.  To be more explicit, we put
  \[
  \tau{[a,b]} = (-b,-a),\quad\tau{(a,b)} = [-b,-a],\quad\tau{[a,b)} =
  [-b,-a),\quad \tau{(a,b]} = (-b,-a].
  \]
  Then the correspondence $J \mapsto \tau{J}$ extends to an involution
  $\tau$ of $I(\R)$
  \[
  J_1 \cup \cdots \cup J_r \mapsto \tau{J_r} \cup \cdots \cup
  \tau{J_1}.
  \]
  Let $\alpha \co  \R \to (0,1)$ be an order preserving
  homeomorphism and let
  \[
  \alpha_t(s) = (1-t)s + t\alpha(s)
  \]
  for $t \in [0,1]$ and $s \in \R$.  Since $\alpha_t \co  \R \to \R$
  is an embedding, it induces a map of partial monoids $I(\alpha_t)
  \co  I(\R) \to I(\R)$ for every $t$, and hence we can define a
  homotopy $h_t \co  I(\R) \to I(\R)^2$ by
  \[
  h_t(P) = (I(\alpha_t)(P),\tau I(\alpha_t)(P)).
  \]
  Clearly, $h_t$ is a map of partial monoids and we have $h_0 =
  (1,\tau)$ because $I(\alpha_0)$ is the identity.
  On the other hand, $h_1$ maps $I(\R)$ into $I(\R)_2$ because
  $I(\alpha)(P)$ is contained in $(0,1)$ and hence is disjoint from
  $\tau I(\alpha)(P) \subset (-1,0)$.
  Finally, we can define a homotopy $\Sigma h_1 \simeq \emptyset$ by
  moving $I(\alpha)(P)$ to negative direction and $\tau I(\alpha)(P)$
  to positive direction, simultaneously, so that the strings $J$ in
  $I(\alpha)(P)$ meet with $\tau J$ at the origin and the resulting
  half-open intervals eventually vanish.
\end{proof}

Let $I(\R)_+$ be the subset of $I(\R)$ consisting of those $J_1 \cup
\cdots \cup J_r$ such that every $J_i$ is a closed interval.  Clearly,
$I(\R)_+$ is a partial submonoid of $I(\R)$.

\begin{Definition}
  Given an orthogonal $G$--module $V$ and a $G$--partial monoid $M$, let
  \[
  I(V,M) = C(V,I(\R) \wedge M), \quad
  I_+(V,M) = C(V,I(\R)_+ \wedge M)\,.
  \]
\end{Definition}


For any $G$--partial monoid $M$, $I(\R) \wedge M$ is a homotopically
invertible $G$--partial monoid with homotopy inversion $\tau \wedge 1$.
Thus we have the following proposition.

\begin{Proposition}
  \label{prop:I(V,M)-is-grouplike}
  If $V$ is a $G$--universe then $I(V,M)$ is grouplike for any $M$.
\end{Proposition}
\section[Proof of Theorem~\ref{thm:group-completion}]{Proof of \fullref{thm:group-completion}}

To establish a relation between $I(V,M)$ and $C(V,M)$, let us choose a
linear embedding $e \co  \R \to V^G$ and a $G$--linear isometry $l
\co  V \times V \to V$.  Then we can define
\[
\lambda \co  I_+(V,M) \to C(V,M)
\]
to be the map which sends a finite set of labeled strings $\{(\{v_i\}
\times J_i,a_i)\}$ to the set of labeled particles
$\{(l(v_i,e(\hat{J}_i)),a_i)\}$, where $\hat{J}_i$ is the middle point
of the closed interval $J_i$.
Note that $\smash{(v_i,e(\hat{J}_i))}$ are pairwise distinct, hence so are
$\smash{l(v_i,e(\hat{J}_i))}$.

\begin{Proposition}
  \label{prop:Iplus_to_C}
  $\lambda \co  I_+(V,M) \to C(V,M)$ is a $G$--homotopy equivalence
  of Hopf $G$--spaces.
\end{Proposition}

\begin{proof}
  Since $\lambda$ is natural with respect to $M$, it extends to a map
  of $G$--equivariant $\Gamma$--spaces.  This, of course, implies that
  $\lambda$ is a map of Hopf $G$--spaces.

  To see that $\lambda$ is a $G$--homotopy equivalence, let
  \(
  \gamma \co  C(V,M) \to I_+(V,M)
  \)
  be a pointed $G$--map which sends a finite set of labeled particles
  $\{(v_i,a_i)\}$ to the set of labeled strings $\{(\{v_i\} \times
  [-1,1],a_i)\}$.  Then we have
  \[
  \gamma\lambda(\{(\{v_i\} \times J_i,a_i)\}) =
  \{(\{l(v_i,e(\hat{J}_i))\} \times [-1,1],a_i)\}
  \]
  and we can define a $G$--homotopy $\gamma\lambda \simeq 1$ by
  \[
  (\gamma\lambda)_t(\{(\{v_i\} \times J_i,a_i)\}) =
  \begin{cases}
    \{(\{l(v_i,e_{2t}(\hat{J}_i))\} \times
    \mathcal{I}_{2t}(J_i),a_i)\},& 0 \leq t \leq 1/2
    \\
    \{(\{l_{2t-1}(v_i)\} \times J_i,a_i)\},& 1/2 \leq t \leq 1
  \end{cases}
  \]
  where
  \begin{enumerate}
  \item $e_t \co  \R \to V^G$ is a linear map $s \mapsto (1-t)e(s)$.
  \item If $J = [a,b]$ then $\mathcal{I}_t(J) = [ta-(1-t),tb+(1-t)]$.
    Thus $\{\mathcal{I}_t(J)\}$ is a continuous family of closed
    intervals such that $\mathcal{I}_0(J) = [-1,1]$ and
    $\mathcal{I}_1(J) = J$.
  \item $\{l_t\}$ is a continuous family of $G$--linear isometries $V
    \to V$ such that $l_0(v) = l(v,0)$ and $l_1$ is the identity of
    $V$.  (Such a family certainly exists because the space of
    $G$--linear isometries $V \to V$ is contractible if $V$ is a
    $G$--universe.)
  \end{enumerate}
  On the other hand, we can define a $G$--homotopy $\lambda\gamma
  \simeq 1$ by
  \[
  (\lambda\gamma)_t(\{(v_i,a_i)\}) = \{(l_t(v_i),a_i)\}.\proved
  \]
\end{proof}

Now let $\rho \co  I_+(V,M) \to I(V,M)$ be the map induced by the
inclusion $I(\R)_+ \subset I(\R)$.  To complete the proof of
\fullref{thm:group-completion}, we need to show that
\begin{equation}
  \label{eq:rho^H}
  \rho^H \co  I_+(V,M)^H \to I(V,M)^H
\end{equation}
is a group completion for every subgroup $H$ of $G$.
Since $V$ is an $H$--universe for any subgroup $H$ of $G$, we need only
consider the case $H = G$.  But then we have:

\begin{Lemma}
  \label{lem:I_+ to I}
  $I_+(V,M)^G \to I(V,M)^G$ is a group completion for a $G$--partial
  monoid $M$ if so is $I_+(\Rinfty,M) \to I(\Rinfty,M)$ for all
  partial abelian monoids $M$.
\end{Lemma}

\begin{proof}
  Let $\mathcal{F}$ be a family of orbit types and let
  $\smash{C(V,M)_{\mathcal{F}}}$ be the subspace of $C(V,M)$ consisting of those
  $(c,v)\in C(V,M)$ such that $\smash{c\in V_{\mathcal{F}}}$, where
  $\smash{V_{\mathcal{F}}=\{v\in V| G_v\in \mathcal{F}\}}$.
  If $\smash{\mathcal{F}_1}$ and $\smash{\mathcal{F}_2}$ are families such
  that $\smash{\mathcal{F}_1\subset \mathcal{F}_2}$ and
  $\smash{\mathcal{F}_2-\mathcal{F}_1}$ consists of only one conjugacy class
  $(H)$
  then we have a fibration sequence
  \[
  C(V,M)_{\mathcal{F}_1}^G\to C(V,M)_{\mathcal{F}_2}^G\to
  C(V,M)_{(H)}^G.
  \]
  Therefore, we see that $\smash{I_+(V,M)^G \to
  I(V,M)^G}$ is a group completion if and only if so are
  $\smash{I_+(V,M)_{(H)}^G\to I(V,M)_{(H)}^G}$, by arguing as in \S{6} of \mbox{Caruso and Waner
  \cite{caruso-waner:approximation}}. But we have
  \[
  C(V,M)_{(H)}^G\cong C(V^H,M^H)_{(H)}^{NH} \simeq C(\Rinfty,
  EW(H)\wedge_{W(H)} M^H ).
  \]
  It follows that $I_+(V,M)^G \to I(V,M)^G$ is a group completion if
  so are
  \[
  I_+(\Rinfty,EW(H)\wedge_{W(H)} M^H) \to I(\Rinfty,EW(H)\wedge_{W(H)}
  M^H).\proved
  \]
\end{proof}

In order to prove the lemma in the non equivariant case we need a
CW--monoid replacement for $I(\Rinfty,M)$.
For any $M$ let $\RS{M}$ be the realization of the total singular
complex of $M$ regarded as a partial abelian monoid such that
\[
\RS{M}_n = \RS{M_n} \subset \RS{M}^n \quad (n \geq 0).
\]
Let $D(M)$ be the classifying space of the permutative category
$\smash{\cQ{\RS{M}}}$ whose space of objects is $\smash{\coprod_{p \geq 0} \RS{M}^p}$
and whose morphisms from $\smash{(a_i) \in \RS{M}^p}$ to $\smash{(b_j) \in \RS{M}^q}$
are given by a map of finite sets $\theta \co  \{1,\dots,p\} \to
\{1,\dots,q\}$ such that $\smash{b_j = \sum_{i \in \theta^{-1}(j)} a_i}$ hold.
Then $D(M)$ is a CW--monoid since it is homeomorphic to the
realization of the diagonal simplicial set $[n] \mapsto N_n\cQ{S_nM}$.
Moreover, there is a natural weak equivalence of Hopf spaces $\Phi
\co  D(M) \to C(\Rinfty,M)$.  (For details, see
Shimakawa \cite[\S{2.4}]{completion}.)
Thus to prove \fullref{lem:I_+ to I} we need only show the following

\begin{Proposition}
  \label{prop:D(IR,M)}
  The natural map $D(I(\R)_+ \wedge M) \to D(I(\R) \wedge M)$ induced
  by the inclusion $I(\R)_+ \subset I(\R)$ is a group completion.
\end{Proposition}

The rest of this section is devoted to the proof of this proposition.

Given a map of topological monoids $f \co  D \to D'$ let $B(D,D')$
denote the realization of the category $\mathcal{B}(D,D')$ whose space
of objects is $D'$ and whose space of morphisms is the product $D
\times D'$, where $(d,d') \in D \times D'$ is regarded as a morphism
from $d'$ to $f(d) \cdot d'$.
Then there is a sequence of maps
\[
D' = B(0,D') \to B(D,D') \to B(D,0) = \BD
\]
induced by the maps $0 \to D$ and $D' \to 0$ respectively.
Observe that $\BD$ is the standard classifying space of the monoid $D$
and $B(D,D)$ is contractible when $f$ is the identity.

In particular, let us take $D = D(I(\R)_+ \wedge M)$ and $D' = D(I(\R)
\wedge M)$, and let $i \co  D \to D'$ be the monoid map induced by
the inclusion $I(\R)_+ \to I(\R)$.  Then there is a commutative
diagram
\begin{equation}
  \label{eq:ladder}
  \begin{CD}
    D @>>> B(D,D) @>>> BD \\
    @V{i}VV @VV{B(1,i)}V @| \\
    D' @>>> B(D,D') @>>> BD
  \end{CD}
\end{equation}
in which the upper and the lower sequences are associated with the
identity and the inclusion $i \co  D \to D'$, respectively.

\begin{Lemma}
  \label{lem:BD}
  The natural map $D \to \Loop{BD}$ is a group completion.
\end{Lemma}

This follows from the fact that $D$ is a homotopy commutative monoid.

\begin{Lemma}
  \label{lem:B(D,D')}
  The lower sequence in the diagram {\eqref{eq:ladder}} is a
  homotopy fibration sequence with contractible total space. 
\end{Lemma}

\propref{prop:D(IR,M)} is deduced from this, because $D \to D'$ is
equivalent to the group completion map $D \to \Loop{BD}$ under the
equivalence $D' \simeq \Loop{BD}$.

\begin{proof}[Proof of \fullref{lem:B(D,D')}]
  By \fullref{prop:I(V,M)-is-grouplike}, $D' = D(I(\R) \wedge
  M)$ is grouplike with homotopy inverse induced by the homotopy
  inversion $\tau \wedge 1$.  Hence $D$ acts on $D'$ through homotopy
  equivalences, and the diagram
  \begin{equation*}
    \label{eq:cartesian}
    \begin{CD}
      D' @>>> B(D,D') \\
      @VVV @VVV \\
      0 @>>> B(D,0)
    \end{CD}
  \end{equation*}
  is homotopy cartesian by Proposition~1.6 of Segal \cite{segal:category}.
  This implies that the lower sequence in the diagram
  \eqref{eq:ladder} is a homotopy fibration sequence.

  It remains to prove that $B(D,D')$ is contractible.  In
  \cite{completion}, we proved this in the case where the partial
  monoid $X \wedge \pm{M}$ is strictly invertible and is generated by
  the elements of $X \wedge M$ and their inverses.
  But the argument given there still applies to the current case, once
  we make the following change in the notation.

  Replace $X \wedge M$ and $X \wedge \pm{M}$ by $I_+(\R) \wedge M$ and
  $I(\R) \wedge M$, respectively, and for any $S = (P_j \wedge a_j)
  \in S_0(I(\R) \wedge M)^p$ put
  \[
  S_+ = (P_j^+ \wedge a_j), \quad S_- = (P_j^- \wedge a_j), \quad
  \overline{S} = (\tau P_j \wedge  a_j),
  \]
  where $\smash{P_j^+}$ and $\smash{P_j^-}$ are the unions of closed intervals and of
  open or half-open intervals contained in $P_j$, respectively.  Note
  that we have $\smash{P_j = P_j^+ \bigcup P_j^-}$ and $\smash{P_j^+ \in I_+(\R)}$.
  Also, for any $S$ such that $S = S_-$ the path
  $[S] \to [\mathbf{0}^p]$ in $B(D,D')$
  should be defined to be the composite
  \[
  [S] \to [\overline{S}_+ \cdot S] \to
  [\overline{I(\alpha)_*(S)}_+ \cdot I(\alpha)_*(S)]
  \xrightarrow{\nabla} [\mathbf{0}^p]
  \]
  where $\alpha$ is a homeomorphism $\R \cong (0,1)$, $I(\alpha)_*(P_j
  \wedge a_j) = (I(\alpha)(P_j) \wedge a_j)$, and $\nabla$ is induced
  by the homotopy
  \[
  \textstyle
  \tau I(\alpha)(P_j)^+ \wedge a_j + I(\alpha)(P_j) \wedge a_j = (\tau
  I(\alpha)(P_j)^+ \bigcup I(\alpha)(P_j)) \wedge a_j \simeq \emptyset
  \wedge a_j = 0.
  \]
  (Compare the proof of
  \fullref{lem:IR-is-homotopically-invertible}.)
\end{proof}
\section[Proof of \ref{thm:infinite-delooping}]{Proof of \fullref{thm:infinite-delooping}}

By a simplicial pointed $G$--space we shall mean a simplicial object in
the category of pointed $G$--spaces and basepoint preserving $G$--maps.
If $\simp{X}$ is a simplicial pointed $G$--space then the basepoints of
$X_n$ form the simplicial set $*$.  Let
\[
\thickr{\simp{X}}' = \thickr{\simp{X}}/\thickr{*}.
\]
Then the natural map $\thickr{\simp{X}} \to \thickr{\simp{X}}'$ is a
$G$--homotopy equivalence, and the maps $\Delta^n \times X_n \to
\thickr{\simp{X}}$ induce \( \Delta^n_+ \wedge X_n \to
\thickr{\simp{X}}'.  \)

Let $T$ be a $G$--equivariant continuous functor $\Top(G) \to \Top(G)$.
Then any simplicial pointed $G$--space $\simp{X}$ is associated with a
$G$--map $\smash{\thickr{T(\simp{X})}' \to T(\thickr{\simp{X}}')}$ induced by
the maps
\[
\Delta^n \times T(X_n) \to \Delta^n_+ \wedge T(X_n) \to
T(\Delta^n_+ \wedge X_n) \to T(\thickr{\simp{X}}')\,.
\]
The next proposition plays a key role in the proof of
\fullref{thm:infinite-delooping}.

\begin{Proposition}
  Let $T \co  \Top(G) \to \Top(G)$ be a $G$--equivariant continuous
  functor. Suppose $T$ satisfies the conditions below.
  \begin{enumerate}
  \item[(C1)] $T(*) = *$.
  \item[(C2)] For any simplicial pointed $G$--space $X_{\bullet}$ the
    natural map $\thickr{T(X_{\bullet})}' \to
    T(\thickr{X_{\bullet}}')$ is a $G$--homotopy equivalence.
  \item[(C3)] For any $X$ and $Y$ the map $T(X \vee Y) \to T(X) \times
    T(Y)$ induced by the projections $X \vee Y \to X$ and $X \vee Y
    \to Y$ is a $G$--homotopy equivalence.
  \item[(C4)] For any subgroup $H$ the natural map $T(G/H_+ \wedge X)
    \to \Map_0(G/H_+,T(X))$, whose adjoint $G/H_+ \wedge T(G/H_+
    \wedge X) \to T(X)$ is induced by the pairing $G/H_+ \wedge G/H_+
    \wedge X \to X$ which sends $(s,s,x)$ to $x$ and $(s,t,x)$ ($s
    \neq t$) to the basepoint of $X$, is a $G$--homotopy equivalence.
  \end{enumerate}
  Suppose further that $T(X)^H$ is grouplike for any $X$ and any
  subgroup $H$ of $G$. Then the following hold.
  \begin{enumerate}
  \item For any orthogonal $G$--module $W$ the natural map $T(X) \to
    \Omega^W T(\Sigma^W X)$ adjoint to $S^W \wedge T(X) \to T(S^W
    \wedge X)$ is a weak $G$--homotopy equivalence.
  \item The correspondence $X \mapsto \left\{ \pi_n T(X)^G \right\}$
    is extendible to an $RO(G)$--graded equivariant homology theory
    defined on the category of pointed $G$--spaces.
  \end{enumerate}
\end{Proposition}

\begin{proof}
  For any pointed $G$--space $X$ let $E(X) = \Omega{T(\Sigma{X})}$.  If
  $T$ satisfies (C1), (C2) and (C3) then by the equivariant version of
  \cite[Theorem~2.12]{config} the natural map $T(X) \to E(X)$ is a
  $G$--equivariant group completion and the sequence
  \[ E(A) \to E(X) \to E(X \cup CA) \]
  associated with a pair of pointed $G$--spaces $(X,A)$ is a
  $G$--fibration sequence up to weak $G$--equivalence.
  But $T(X) \to E(X) = \Omega{T(\Sigma{X})}$ is a weak $G$--equivalence
  because $T(X)^H$ is grouplike for any subgroup $H$.  Hence
  \[ T(A) \to T(X) \to T(X \cup CA) \]
  is a $G$--fibration sequence up to weak $G$--equivalence.
  Moreover, $T$ preserves $G$--homotopies because it is a
  $G$--equivariant continuous functor.
  Therefore, the correspondence $X \mapsto \left\{ \pi_n T(X)^G
  \right\}$ determines a $\mathbb{Z}$--graded equivariant homology
  theory.

  Let $\Gamma_G$ be the full subcategory of $\Top(G)$ consisting of
  finite pointed $G$--sets. To prove the assertions we need only show
  that the correspondence $S \mapsto T(S \wedge X)$ from $\Gamma_G$ to
  $\Top(G)$ is a special $\Gamma_G$--space in the sense of
  \cite{infinite}.  But this follows from the conditions (C3) and
  (C4).
\end{proof}

Now let $T(X) = I(V,X \wedge M)$.  We shall show that $T$ satisfies
the conditions (C1), (C2), (C3) and (C4).  This of course proves
\fullref{thm:infinite-delooping}.

It is obvious that (C1) holds.  (C2) is proved by the argument similar
to the one used in the proof of \cite[Theorem 3.2]{config}.
To prove (C3) let us define
\[
T(X) \times T(Y) \to T(X \vee Y)
\]
to be the composite
\[
I(V,X \wedge M) \times I(V,Y \wedge M) \xrightarrow{(i_*,j_*)}%
I(V,(X \vee Y) \wedge M)^2 \xrightarrow{\mu} I(V,(X \vee Y) \wedge M)
\]
where $i_*$ and $j_*$ are induced by the inclusions of $X$ and $Y$
into $X \vee Y$, respectively, and $\mu$ is the multiplication of the
Hopf $G$--space $I(V,(X \vee Y) \wedge M)$.
By using the fact that the space of $G$--linear isometries of $V$ is
contractible one can show that the map above gives a $G$--homotopy
inverse to $T(X \vee Y) \to T(X) \times T(Y)$.
Finally, to prove (C4) let us choose a $G$--embedding $G/H \to V$ and a
$G$--linear isometry $l \co  V \times V \to V$.  Then we can
construct a $G$--homotopy inverse to the natural map $T(G/H_+ \wedge X)
\to \Map_0(G/H_+,T(X))$ by the following procedure:

\begin{enumerate}
  
\item For given $f \co  G/H_+ \to T(X)$ let us write
  \[
  f(gH) = (c(gH),P(gH) \wedge a(gH))
  \]
  where $c(gH) \subset V$, $P(gH) \co  c(gH) \to I(\R)$ and $a(gH)
  \co  c(gH) \to X \wedge M$.
\item Let $\tilde{c}$ be the image of the union $\bigcup \{gH\} \times
  c(gH)$ under the embedding
  \[
  \iota \co  G/H \times V \subset V \times V \xrightarrow{l} V.
  \]
\item Define $\tilde{a} \co  \tilde{c} \to I(\R) \wedge G/H_+ \wedge
  X \wedge M$ by
  \[
  \tilde{a}(\iota(gH,\xi)) = P(gH) \wedge gH \wedge a(gH)(\xi), \quad
  \xi \in c(gH).
  \]
\item Define $\rho \co  \mbox{Map}(G/H,T(X)) \to T(G/H_+ \wedge X)$
  by $\rho(f) = (\tilde{c},\tilde{a})$.
\end{enumerate}
That $\rho$ gives a $G$--homotopy inverse to $T(G/H_+ \wedge X) \to
\Map(G/H,T(X))$ follows, again, from the contractibility of the space
of $G$--linear isometries of $V$.

\bibliographystyle{gtart}
\bibliography{link}

\end{document}